\documentclass{conm-p-l} 

\usepackage{amsmath,amsfonts,amssymb,amsthm,amscd,epsfig,comment}
\usepackage[all]{xy}

\def\C{\mathbb{C}}

\def\Z{\mathbb{Z}}

\def\R{{\mathbb{R}}}

\def\g{\ensuremath{\mathfrak{g}}}

\def\v{\mathbf{v}}
\def\w{\mathbf{w}}
\def\e{\mathfrak{e}}
\def\E{\mathcal{E}}

\DeclareMathOperator{\im}{im} 
\DeclareMathOperator{\Hom}{Hom}

\DeclareMathOperator{\inc}{in}
\DeclareMathOperator{\out}{out}

\DeclareMathOperator{\Irr}{Irr}

\DeclareMathOperator{\wt}{wt}
\DeclareMathOperator{\Span}{Span}

\DeclareMathOperator{\rep}{rep}
\DeclareMathOperator{\md}{mod}
\DeclareMathOperator{\res}{res}

\newtheorem{theo}{Theorem}[section]
\newtheorem{prop}[theo]{Proposition}

\newtheorem{cor}[theo]{Corollary}

\newtheorem{defin}[theo]{Definition}
\newtheorem*{rem*}{Remark}
\newtheorem{rem}[theo]{Remark}

\numberwithin{equation}{section}

\begin{document}
\title{Quivers and the Euclidean group}
\author{Alistair Savage}
\address{University of Ottawa\\
Ottawa, Ontario \\ Canada} \email{alistair.savage@uottawa.ca}
\thanks{This research was supported in part by the Natural
Sciences and Engineering Research Council (NSERC) of Canada}
\date{December 10, 2007}

\begin{abstract}
We show that the category of representations of the Euclidean group
of orientation-preserving isometries of two-dimensional Euclidean
space is equivalent to the category of representations of the
preprojective algebra of type $A_\infty$.  We also consider the
moduli space of representations of the Euclidean group along with a
set of generators.  We show that these moduli spaces are quiver
varieties of the type considered by Nakajima.  Using these
identifications, we prove various results about the representation
theory of the Euclidean group.  In particular, we prove it is of
wild representation type but that if we impose certain restrictions
on weight decompositions, we obtain only a finite number of
indecomposable representations.
\end{abstract}

\maketitle

\section{Introduction}

The Euclidean group $E(n) = \R^n \rtimes SO(n)$ is the group of
orientation-preserving isometries of $n$-dimensional Euclidean
space. The study of these objects, at least for $n=2,3$, predates
even the concept of a group. In this paper we will focus on the
Euclidean group $E(2)$. Even in this case, much is still unknown
about the representation theory.

Since $E(2)$ is solvable, all its finite-dimensional irreducible
representations are one-dimensional.  The finite-dimensional unitary
representations, which are of interest in quantum mechanics, are
completely reducible and thus isomorphic to direct sums of such
one-dimensional representations. The infinite-dimensional unitary
irreducible representations have received considerable attention
(see \cite{AD00,BAU94,BS00}). There also exist finite-dimensional
nonunitary indecomposable representations (which are not
irreducible) and much less is known about these. However, they play
an important role in mathematical physics and the representation
theory of the Poincar\'e group.  The Poincar\'e group is the group
of isometries of Minkowski spacetime. It is the semidirect product
of the translations of $\R^3$ and the Lorentz transformations.  In
1939, Wigner \cite{Wig39} studied the subgroups of the Lorentz group
leaving invariant the four-momentum of a given free particle.  The
maximal such subgroup is called the \emph{little group}.  The little
group governs the internal space-time symmetries of the relativistic
particle in question.  The little groups of massive particles are
locally isomorphic to the group $O(3)$ while the little groups of
massless particles are locally isomorphic to $E(2)$.  That is, their
Lie algebras are isomorphic to those of $O(3)$ and $E(2)$
respectively.  We refer the reader to \cite{Bar86,COT89,KN86,Mir94}
for further details.

The group $E(2)$ also appears in the Chern-Simons formulation of
Einstein gravity in $2+1$ dimensions.  In the case when the
space-time has Euclidean signature and the cosmological constant
vanishes, the phase space of gravity is the moduli space of flat
$E(2)$-connections.

In the current paper, we relate the representation theory of the
Euclidean group $E(2)$ to the representation theory of preprojective
algebras of quivers of type $A_\infty$.  In fact, we show that the
categories of representations of the two are equivalent.  To prove
this, we introduce a modified enveloping algebra of the Lie algebra
of $E(2)$ and show that it is isomorphic to the preprojective
algebra of type $A_\infty$.  Furthermore, we consider the moduli
space of representations of $E(2)$ along with a set of generators.
We show that these moduli spaces are quiver varieties of the type
considered by Nakajima in \cite{Nak94,Nak98}. These identifications
allow us to draw on known results about preprojective algebras and
quiver varieties to prove various statements about representations
of $E(2)$.  In particular, we show that $E(2)$ is of wild
representation type but that if we impose certain restrictions on
the weight decomposition of a representation, we obtain only a
finite number of indecomposable representations.  We conclude with
some potential directions for future investigation.


\section{The Euclidean algebra}

Let $E(2) = \R^2 \rtimes SO(2)$ be the Euclidean group of motions in
the plane and let $\e(2)$ be the complexification of its Lie
algebra.  We call $\e(2)$ the (three-dimensional) Euclidean algebra.
It has basis $\{p_+, p_-, l\}$ and commutation relations
\begin{equation} \label{eq:e2-relations}
[p_+, p_-]=0,\quad [l,p_\pm] = \pm p_\pm.
\end{equation}
Since $SO(2)$ is compact, the category of finite-dimensional
$E(2)$-modules is equivalent to the category of finite-dimensional
$\e(2)$-modules in which $l$ acts semisimply with integer
eigenvalues.  Will will use the term $\e(2)$-module to refer only to
such modules.  For $k \in \Z$, we shall write $V_k$ to indicate the
eigenspace of $l$ with eigenvalue $k$ (the \emph{$k$-weight space}).
Thus, for an $\e(2)$-module $V$, we have the weight space
decomposition
\[
V= \bigoplus_k V_k,\quad V_k = \{v \in V\ |\ l \cdot v = kv\},\quad
k \in \Z,
\]
and
\[
p_+ V_k \subseteq V_{k+1},\quad p_- V_k \subseteq V_{k-1}.
\]

We may form the tensor product of any representation $V$ with the
character $\chi_n$ for $n \in \Z$.  Here $\chi_n$ is the
one-dimensional module $\C$ on which $p_\pm$ act by zero and $l$
acts by multiplication by $n$.  Then a weight space $V_k$ of weight
$k$ becomes a weight space $V_k \otimes \chi_n$ of weight $k + n$.
In this way, we may ``shift weights'' as we please.

For $k \in \Z$, let $\mathbf{e}^k$ be the element of $(\Z_{\ge
0})^\Z$ with $k$th component equal to one and all others equal to
zero.  For an $\e(2)$-module $V$ we define
\[
\mathbf{dim}\, V = \sum_{k \in \Z} (\dim V_k) \mathbf{e}^k.
\]

Let $U$ be the universal enveloping algebra of $\e(2)$ and let
$U^+$, $U^-$ and $U^0$ be the subalgebras generated by $p_+$, $p_-$
and $l$ respectively.  Then we have the triangular decomposition
\[
U \cong U^+ \otimes U^0 \otimes U^- \quad \text{(as vector spaces)}.
\]
Note that the category of representations of $U$ is equivalent to
the category of representations of $\e(2)$.  In
\cite[Chapter~23]{L93}, Lusztig introduced the modified quantized
enveloping algebra of a Kac-Moody algebra.  Following this idea, we
introduce the \emph{modified enveloping algebra} $\tilde U$ by
replacing $U^0$ with a sum of 1-dimensional algebras
\[
\tilde U = U^+ \otimes \left( \bigoplus_{k \in \Z} \C a_k \right)
\otimes U^-.
\]
Multiplication is given by
\begin{gather*}
a_k a_l = \delta_{kl} a_k, \\
p_+ a_k = a_{k+1} p_+,\quad p_- a_k = a_{k-1} p_-, \\
p_+ p_- a_k = p_- p_+ a_k.
\end{gather*}
One can think of $a_k$ as projection onto the $k$th weight space.
Note that $\tilde U$ is an algebra without unit.  We say a $\tilde
U$-module $V$ is \emph{unital} if
\begin{enumerate}
\item for any $v \in V$, we have $a_k v = 0$ for all but finitely
many $k \in Z$, and
\item for any $v \in V$, we have $\sum_{k \in \Z} a_k v = v$.
\end{enumerate}
A unital $\tilde U$-module can be thought of as a $U$-module with
weight decomposition.  Thus we have the following proposition.

\begin{prop} \label{prop:tU=U}
The category of unital $\tilde U$-modules is equivalent to the
category of $U$-modules and hence the category of $\e(2)$-modules.
\end{prop}


\section{Preprojective algebras}
\label{sec:preproj}

In this section, we review some basic results about preprojective
algebras.  The reader is referred to \cite{GLS04} for further
details.

A \emph{quiver} is a 4-tuple $(I,H, \out, \inc)$ where $I$ and $H$
are disjoint sets and $\out$ and $\inc$ are functions from $H$ to
$I$. The sets $I$ and $H$ are called the \emph{vertex set} and
\emph{arrow set} respectively.  We think of an element $h \in H$ as
an arrow from the vertex $\out(h)$ to the vertex $\inc(h)$.
\[
\xymatrix{ *=0{\bullet} \ar[rr]^<{\out(h)}^>{\inc(h)}^h & &
*=0{\bullet}}
\]
An arrow $h \in H$ is called a \emph{loop} if $\out(h) = \inc(h)$. A
quiver is said to be \emph{finite} if both its vertex and arrow sets
are finite.

We shall be especially concerned with the following quivers.  For
$a, b \in \Z$ with $a \le b$, let $Q_{a,b}$ be the quiver with
vertex set $I=\{k \in \Z\ |\ a \le k \le b\}$ and arrows $H = \{h_i\
|\ a \le i \le b-1\}$ with $\out(h_i) = i$ and $\inc(h_i) = i+1$. We
say that $Q_{a,b}$ is a quiver of type $A_{b-a+1}$ since this is the
type of its underlying graph.  The quiver $Q_\infty$ has vertex set
$I=\Z$ and arrows $H = \{h_i\ |\ i \in \Z\}$ with $\out(h_i) = i$
and $\inc(h_i) = i+1$.  We say that the quiver $Q_\infty$ is of type
$A_\infty$.  Note that the quivers $Q_{a,b}$ are finite while the
quiver $Q_\infty$ is not.

Let $Q=(I,H,\out,\inc)$ be a quiver without loops and let
$Q^*=(I,H^*,\out^*,\inc^*)$ be the \emph{double quiver} of $Q$. By
definition,
\begin{gather*}
H^* = \{h\ |\ h \in H\} \cup \{\bar h\ |\ h \in H\},\\
\out^*(h) = \out(h),\quad \inc^*(h) = \inc(h),\quad \out^*(\bar h) =
\inc(h),\quad \inc^*(\bar h) = \out(h).
\end{gather*}
From now on, we will write $\inc$ and $\out$ for $\inc^*$ and
$\out^*$ respectively.  Since $\inc^*|_H = \inc$ and $\out^*|_H =
\out$, this should cause no confusion.

A \emph{path} in a quiver $Q$ is a sequence $p = h_n h_{n-1} \cdots
h_1$ of arrows such that $\inc(h_i) = \out(h_{i+1})$ for $1 \le i
\le n-1$.  We call the integer $n$ the \emph{length} of $p$ and
define $\out(p) = \out(h_1)$ and $\inc(p) = \inc(h_n)$. The
\emph{path algebra} $\C Q$ is the algebra spanned by the paths in
$Q$ with multiplication given by
\[
p \cdot p' = \begin{cases} pp' & \text{if } \inc(p') = \out(p) \\
0 & \text{otherwise} \end{cases}
\]
and extended by linearity.  We note that there is a trivial path
$\epsilon_i$ starting and ending at $i$ for each $i \in I$.  The
path algebra $\C Q$ has a unit (namely $\sum_{i \in I} \epsilon_i$)
if and only if the quiver $Q$ is finite.

A \emph{relation} in a quiver $Q$ is a sum of the form $\sum_{j=1}^k
a_j p_j$, $a_j \in \C$, $p_j$ a path for $1 \le j \le k$. For $i \in
I$ let
\[
r_i = \sum_{h \in H,\, \out(h)=i} \bar h h - \sum_{h \in H,\,
\inc(h) = i} h \bar h
\]
be the \emph{Gelfand-Ponomarev relation} in $Q^*$ associated to $i$.
The \emph{preprojective algebra} $P(Q)$ corresponding to $Q$ is
defined to be
\[
P(Q) = \C Q^*/J
\]
where $J$ is the two-sided ideal generated by the relations $r_i$
for $i \in I$.

Let $\mathcal{V}(I)$ denote the category of finite-dimensional
$I$-graded vector spaces with morphisms being linear maps respecting
the grading.  For $V \in \mathcal{V}(I)$, we let $\mathbf{dim}\, V =
(\dim V_i)_{i \in I}$ be the $I$-graded dimension of $V$.  A
\emph{representation} of the quiver $Q^*$ is an element $V \in
\mathcal{V}(I)$ along with a linear map $x_h : V_{\out(h)} \to
V_{\inc(h)}$ for each $h \in H^*$. We let
\[
\rep (Q^*,V) = \bigoplus_{h \in H^*} \Hom_\C (V_{\out(h)},
V_{\inc(h)})
\]
be the affine variety consisting of representations of $Q^*$ with
underlying graded vector space $V$.  A representation of a quiver
can be naturally interpreted as a $\C Q^*$-module structure on $V$.
For a path $p = h_n h_{n-1} \dots h_1$ in $Q^*$, we let
\[
x_p = x_{h_n} x_{h_{n-1}} \cdots x_{h_1}.
\]
We say a representation $x \in \rep(Q^*,V)$ \emph{satisfies the
relation} $\sum_{j=1}^k a_j p_j$, if
\[
\sum_{j=1}^k a_j x_{p_j}=0.
\]
If $R$ is a set of relations, we denote by $\rep (Q^*,R,V)$ the set
of all representations in $\rep (Q^*,V)$ satisfying all relations in
$R$. This is a closed subvariety of $\rep(Q^*,V)$.  Every element of
$\rep(Q^*,J,V)$ can be naturally interpreted as a $P(Q)$-module
structure on $V$ and so we also write
\[
\md(P(Q),V) = \rep(Q^*,J,V)
\]
for the affine variety of $P(Q)$-modules with underlying vector
space $V$.

The algebraic group $G_V = \prod_{i \in I} GL(V_i)$ acts on $\md
(P(Q),V)$ by
\[
g \cdot x = (g_i)_{i \in I} \cdot (x_h)_{h \in H^*} = (g_{\inc(h)}
x_h g_{\out(h)}^{-1})_{h \in H^*}.
\]
Two $P(Q)$-modules are isomorphic if and only if they lie in the
same orbit.  For a dimension vector $\v \in (\Z_{\ge 0})^I$, let
\[
V^\v = \bigoplus_{i \in I} \C^{\v_i},\quad \md(P(Q),\v) =
\md(P(Q),V^\v),\quad G_\v = G_{V^\v}.
\]
Then we have that $\md (P(Q),V) \cong \md(P(Q),\mathbf{dim}\, V)$
for all $V \in \mathcal{V}(I)$.  Therefore, we will blur the
distinction between $\md (P(Q),V)$ and  $\md (P(Q),\mathbf{dim}\,
V)$.

We say an element $x \in \md(P(Q),V)$ is \emph{nilpotent} if there
exists an $N \in \Z_{>0}$ such that for any path $p$ of length
greater than $N$, we have $x_p=0$.  Denote the closed subset of
nilpotent elements of $\md(P(Q),V)$ by $\Lambda_{V,Q}$ and let
$\Lambda_{\v,Q} = \Lambda_{V^\v,Q}$. The varieties $\Lambda_{V,Q}$
are called \emph{nilpotent varieties} or \emph{Lusztig quiver
varieties}. Lusztig \cite[Theorem~12.3]{L91} has shown that the
$\Lambda_{V,Q}$ have pure dimesion $\dim (\rep(Q,V))$.

\begin{prop}
For a quiver $Q$, the following are equivalent:
\begin{enumerate}
\item \label{prop-item:PQ-fd} $P(Q)$ is finite-dimensional,
\item \label{prop-item:PQ-nilpotent} $\Lambda_{V,Q} = \md(P(Q),V)$ for all $V
\in \mathcal{V}(I)$,
\item \label{prop-item:Q-dynkin} $Q$ is a Dynkin quiver (i.e. its underlying graph is of $ADE$
type).
\end{enumerate}
\end{prop}

\begin{proof}
The equivalence of \eqref{prop-item:PQ-fd} and
\eqref{prop-item:Q-dynkin} is well-known (see for example
\cite{Rei97}).  That \eqref{prop-item:PQ-nilpotent} implies
\eqref{prop-item:Q-dynkin} was proven by Crawley-Boevey \cite{Cra01}
and the converse was proven by Lusztig \cite[Proposition~14.2]{L91}.
\end{proof}

Thus, for a Dynkin quiver $Q$, nilpotency holds automatically and
$\Lambda_{V,Q}$ is just the variety of representations of the
preprojective algebra $P(Q)$ with underlying vector space $V$.

The representation type of the preprojective algebras is known.
\begin{prop}[{\cite{DR92, GS03}}] \label{prop:preproj-rep-type}
Let $Q$ be a finite quiver.  Then the following hold:
\begin{enumerate}
\item $P(Q)$ is of finite representation type if and only if $Q$ is
of Dynkin type $A_n$, $n \le 4$, and
\item $P(Q)$ is of tame representation type if and only if $Q$ is of
Dynkin type $A_5$ or $D_4$.
\end{enumerate}
Thus $P(Q)$ is of wild representation type if $Q$ is not of Dynkin
type $A_n$, $n \le 5$, or $D_4$.
\end{prop}

In the sequel, we will refer to the preprojective algebra
$P(Q_\infty)$. While $Q_\infty$ is not a finite quiver, any
finite-dimensional representation is supported on finitely many
vertices and thus is a representation of a quiver of type $A_n$ for
sufficiently large $n$. Thus we deduce the following.

\begin{cor} \label{cor:Qinf-type}
All finite-dimensional representations of $P(Q_\infty)$ are
nilpotent and $P(Q_\infty)$ is of wild representation type.
\end{cor}

For a finite quiver $Q$, let $\g_Q$ denote the Kac-Moody algebra
whose Dynkin graph is the underlying graph of $Q$ and let
$U(\g_Q)^-$ denote the lower half of its universal enveloping
algebra. It turns out that Lusztig quiver varieties are intimately
related to $U(\g_Q)^-$.  Namely, Lusztig \cite{L91} has shown that
there is a space of constructible functions on the varieties
$\Lambda_{\v,Q}$, $\v \in (\Z_{\ge 0})^I$, and a natural convolution
product such that this space of functions is isomorphic as an
algebra to $U(\g_Q)^-$. The functions on an individual
$\Lambda_{\v,Q}$ correspond to the weight space of weight $-\sum_{i
\in I} \v_i \alpha_i$, where the $\alpha_i$ are the simple roots of
$\g_Q$.  Furthermore, the irreducible components of $\Lambda_{\v,Q}$
are in one-to-one correspondence with a basis of this weight space.
Under this correspondence, each irreducible component is associated
to the unique function equal to one on an open dense subset of that
component and equal to zero on an open dense subset of all other
components.  The set of these functions yields a basis of
$U(\g_Q)^-$, called the \emph{semicanonical basis}, with very nice
integrality and positivity properties (see \cite{Lus00}).  If
instead of constructible functions one works with the Grothendieck
group of a certain class of perverse sheaves, a similar construction
yields a realization of (the lower half of) the quantum group
$U_q(\g_Q)^-$ and the \emph{canonical basis} (see \cite{L91}).


\section{Representations of the Euclidean algebra and preprojective
algebras}
\label{sec:euc-preproj}

In this section we examine the close relationship between
representations of the Euclidean algebra $\e(2)$ and the
preprojective algebras of types $A_n$ and $A_\infty$.

\begin{theo} \label{thm:tU=PQ}
The modified universal enveloping algebra $\tilde U$ is isomorphic
to the preprojective algebra $P(Q_\infty)$.
\end{theo}

\begin{proof}
Define a map $\psi : \C Q_\infty^* \to \tilde U$ by
\[
\psi(\epsilon_i) = a_i,\quad \psi(h_i) = p_+ a_i = a_{i+1} p_+,\quad
\psi(\bar h_i) = a_i p_- = p_- a_{i+1},\quad i \in I.
\]
It is easily verified that this extends to a surjective map of
algebras with kernel $J$ and thus the result follows.
\end{proof}

Let $\mathbf{Mod}\, \e(2)$ be the category of $\e(2)$-modules.  For
$a \le b$, let $\mathbf{Mod}_{a,b} \, \e(2)$ be the full subcategory
consisting of representations $V$ such that $V_k = 0$ for $k<a$ or
$k> b$.  For $\v \in (\Z_{\ge 0})^\Z$, we also define
$\mathbf{Mod}_{a,b}^\v \, \e(2)$ and $\mathbf{Mod}^\v\, \e(2)$ to be
the full subcategories of $\mathbf{Mod}_{a,b} \, \e(2)$ and
$\mathbf{Mod}\, \e(2)$ consisting of representations $V$ such that
$\mathbf{dim}\, V = \v$.

Let $\mathbf{Mod}\, P(Q)$ be the category of finite-dimensional
$P(Q)$-modules and for $\v \in (\Z_{\ge 0})^I$, let
$\mathbf{Mod}^\v\, P(Q)$ be the full subcategory consisting of
modules of graded dimension $\v$.

\begin{cor}
\label{cor:eucrep=preproj} We have the following equivalences of
categories.
\begin{enumerate}
\item \label{cor-item:eucrep=qmod} $\mathbf{Mod}^\v\, \e(2) \cong \mathbf{Mod}^\v\,
P(Q_\infty)$, $\mathbf{Mod}\, \e(2) \cong \mathbf{Mod}\,
P(Q_\infty)$,
\item \label{cor-item:eucrep=qmod-ab} $\mathbf{Mod}_{a,b}^\v \, \e(2)
\cong \mathbf{Mod}^\v\, P(Q_{a,b})$, $\mathbf{Mod}_{a,b} \, \e(2)
\cong \mathbf{Mod}\, P(Q_{a,b})$.
\end{enumerate}
\end{cor}

\begin{proof}
Statement~\eqref{cor-item:eucrep=qmod} follows from
Theorem~\ref{thm:tU=PQ} and Proposition~\ref{prop:tU=U}.
Statement~\eqref{cor-item:eucrep=qmod-ab} is obtained by restricting
weights to lie between $a$ and $b$.
\end{proof}

\begin{theo} \label{thm:Euc-rep-type}
The following statements hold.
\begin{enumerate}
\item The Euclidean algebra $\e(2)$, and hence the Euclidean group $E(2)$,
have wild representation type, and
\item for $a, b \in \Z$ with $0 \le b-a \le 3$, there are a finite
number of isomorphism classes of indecomposable $\e(2)$-modules $V$
whose weights lie between $a$ and $b$; that is, such that $V_k = 0$
for $k<a$ or $k>b$.
\end{enumerate}
\end{theo}

\begin{proof}
These statements follow immediately from
Corollary~\ref{cor:eucrep=preproj},
Proposition~\ref{prop:preproj-rep-type} and
Corollary~\ref{cor:Qinf-type}.
\end{proof}

\begin{cor} \label{cor:Eucrep-weight-set}
Let $A$ be a finite subset of $\Z$ with the property that $A$ does
not contain any five consecutive integers.  Then there are a finite
number of isomorphism classes of indecomposable $\e(2)$-modules $V$
with the property that $V_k = 0$ if $k \not \in A$.
\end{cor}

\begin{proof}
Partition $A$ into subsets $A_1, \dots, A_n$ such that $A_j = \{a_j,
a_j + 1, \dots, a_j + m_j\}$ and $|a-b| > 1$ for $a \in A_i$, $b \in
A_j$ with $i \ne j$.  By hypothesis, we have $m_j \le 3$ for $1 \le
j \le n$.  Let $V$ be an $\e(2)$-module such that $V_k=0$ if $k \not
\in A$.  Then $V$ decomposes as a direct sum of modules $V =
\bigoplus_{j=1}^n V^j$ where $V^j_k = 0$ if $k<a_j$ or $k > a_j +
m_k$.  Thus, if $V$ is indecomposable, we must have $V = V^j$ for
some $j$.  But there are a finite number of such $V^j$, up to
isomorphism, by Theorem~\ref{thm:Euc-rep-type}.  The result follows.
\end{proof}

For $a \in \Z$, we say an $\e(2)$-module $V$ has \emph{lowest
weight} $a$ if $V_a \ne 0$ and $V_k = 0$ for $k < a$.

\begin{cor} \label{cor:finite-dim5}
For all $a \in \Z$, there are a finite number of isomorphism classes
of indecomposable $\e(2)$-modules with lowest weight $a$ and
dimension less than or equal to five.
\end{cor}

\begin{proof}
By tensoring with the character $\chi_{-a}$ we may assume that
$a=0$. In order for an $\e(2)$-module to be indecomposable, its set
of weights must be a set of consecutive integers.  By
Corollary~\ref{cor:Eucrep-weight-set}, it suffices to consider the
modules of dimension 5.  Again, by
Corollary~\ref{cor:Eucrep-weight-set}, we need only consider the
case when $\dim V_k = 1$ for $0 \le k \le 4$.  We consider the
equivalent problem of classifying the $G_V$-orbits of indecomposable
elements $x \in \Lambda_{V,Q_{0,4}}$ where $V_k = \C$ for $0 \le k
\le 4$. Fixing the standard basis in each $V_k$, we can view the
maps $x_h$, $h \in H^*$, as complex numbers. Considering the
Gelfand-Ponomarev relation $r_0$, we see that $x_{\bar h_0}
x_{h_0}=0$.  Then the relation $r_1$ implies $x_{\bar h_1} x_{h_1} =
0$.  Continuing in this manner, we see that $x_{\bar h_i} x_{h_i} =
0$ for $0 \le i \le 3$.  Thus $x_{h_i} = 0$ or $x_{\bar h_i} = 0$
for $0 \le i \le 3$. Since $x$ is indecomposable, we cannot have
both $x_{h_i}=0$ and $x_{\bar h_i} =0$ for any $i$. Thus, there are
precisely $2^4=16$ $G_V$-orbits in $\Lambda_{V,Q_{0,4}}$.
Representatives for these orbits correspond to setting one of
$x_{h_i}$ or $x_{\bar h_i}$ equal to one and the other to zero for
each $0 \le i \le 3$.
\end{proof}

We note that Douglas \cite{Dou06b} has shown that there are finitely
many indecomposable $\e(2)$-modules (up to isomorphism) of
dimensions five and six.  The proof of
Corollary~\ref{cor:finite-dim5} shows how
Corollary~\ref{cor:Eucrep-weight-set} can simply such proofs.  We
also point out that the graphs appearing in \cite{Dou06b} roughly
correspond, under the equivalence of categories in
Corollary~\ref{cor:eucrep=preproj}, to the diagrams appearing in the
enumeration of irreducible components of quiver varieties given in
\cite{FS03}.

\begin{rem} \label{rem:LQV-rep-e2}
As noted at the end of Section~\ref{sec:preproj}, the Lusztig quiver
varieties $\Lambda_{\v,Q}$ are closely related to the Kac-Moody
algebra $\g_Q$.  Thus, the results of this section show that there
is a relationship between the representation theory of the Euclidean
group $E(2)$ and the Lie algebra $\mathfrak{sl}_\infty$ (or the Lie
groups $SL(n)$).
\end{rem}


\section{Nakajima quiver varieties}
\label{sec:nak-qv}

In this section we briefly review the quiver varieties introduced by
Nakajima \cite{Nak94,Nak98}.  We restrict our attention to the case
when the quiver involved is of type $A$.

Let $Q$ be the quiver $Q_\infty$ or $Q_{a,b}$ for some $a \le b$.
For $V, W \in \mathcal{V}(I)$ define
\[
L_Q(V,W) = \Lambda_{V,Q} \oplus \bigoplus_{i \in I} \Hom_\C (W_i,
V_i).
\]
We denote points of $L_Q(V,W)$ by $(x,s)$ where $x = (x_h)_{h \in
H^*} \in \Lambda_{V,Q}$ and $s = (s_i)_{i \in I} \in \bigoplus_{i
\in I} \Hom_\C (W_i, V_i)$.  We say an $I$-graded subspace $S$ of
$V$ is \emph{$x$-invariant} if $x_h(S_{\out(h)}) \subseteq
S_{\inc(h)}$ for all $h \in H^*$. We say a point $(x,s) \in
L_Q(V,W)$ is \emph{stable} if the following property holds: If $S$
is an $I$-graded $x$-invariant subspace of $V$ containing $\im s$,
then $S=V$.  We denote by $L_Q(V,W)^{\text{st}}$ the set of stable
points.

The group $G_V$ acts on $L_Q(V,W)$ by
\[
g \cdot (x,s) = (g_i)_{i \in I} \cdot ((x_h)_{h \in H^*}, (s_i)_{i
\in I}) = ((g_{\inc(h)} x_h g_{\out(h)}^{-1})_{h \in H^*}, (g_i
s_i)_{i \in I}).
\]
The action of $G_V$ preserves the stability condition and the
stabilizer in $G_V$ of a stable point is trivial.  We form the
quotient
\[
\mathcal{L}_Q(V,W) = L_Q(V,W)^{\text{st}} / G_V.
\]
The $\mathcal{L}_Q(V,W)$ are called \emph{Nakajima quiver
varieties}.  For $\v, \w \in (\Z_{\ge 0})^I$, we set
\[
L_Q(\v,\w) = L_Q(V^\v,V^\w), \quad L_Q(\v,\w)^{\text{st}} =
L_Q(V^\v,V^\w)^{\text{st}}, \quad \mathcal{L}_Q(\v,\w) =
\mathcal{L}_Q(V^\v,V^\w).
\]
We then have
\begin{gather*}
L_Q(V,W) \cong L_Q(\mathbf{dim}\, V, \mathbf{dim}\, W),\quad
L_Q(V,W)^{\text{st}}
\cong L_Q(\mathbf{dim}\, V, \mathbf{dim}\, W)^{\text{st}},\\
\quad \mathcal{L}_Q(V,W) \cong \mathcal{L}_Q(\mathbf{dim}\, V,
\mathbf{dim}\, W),
\end{gather*}
and so we often blur the distinction between these pairs of
isomorphic varieties.

Let $\Irr \Lambda_{V,Q}$ (resp. $\Irr \mathcal{L}_Q(V,W)$) denote
the set of irreducible components of $\Lambda_{V,Q}$ (resp.
$\mathcal{L}_Q(V,W)$).  Then $\Irr \mathcal{L}_Q(V,W)$ can be
identified with
\[
\left\{Y \in \Irr \Lambda_{V,Q}\ \left| \ \left( Y \oplus
\bigoplus_{i \in I} \Hom_\C(W_i,V_i) \right)^{\text{st}} \ne
\emptyset\right.\right\}.
\]
Specifically, the irreducible components of $\mathcal{L}_Q(V,W)$ are
precisely those
\[
\left( \left( Y \oplus \bigoplus_{i \in I} \Hom_\C(W_i,V_i)
\right)^{\text{st}} \right)/G_V
\]
which are nonempty.

\begin{prop}[{\cite[Corollary~3.12]{Nak98}}] \label{prop:NQV-dim}
The dimension of the Nakajima quiver varieties associated to the
quiver $Q_\infty$ are given by
\[
\dim_\C \mathcal{L}_{Q_\infty}(\v,\w) = \sum_{i \in \Z} (\v_i \w_i -
\v_i^2 + \v_i \v_{i+1}).
\]
\end{prop}

In a manner analogous to the way in which Lusztig quiver varieties
are related to $U(\g_Q)^-$ (see Section~\ref{sec:preproj}), Nakajima
quiver varieties are closely related to the representation theory of
$\g_Q$.  In particular, Nakajima \cite{Nak98} has shown that
$\bigoplus_{\v} H_{\text{top}}(\mathcal{L}_Q(\v,\w))$ is isomorphic
to the irreducible integrable highest-weight representation of
$\g_Q$ of highest weight $\sum_{i \in I} \w_i \omega_i$ where the
$\omega_i$ are the fundamental weights of $\g_Q$.  Here
$H_{\text{top}}$ is top-dimensional Borel-Moore homology.  The
action of the Chevalley generators of $\g_Q$ are given by certain
convolution operations. The vector space
$H_{\text{top}}(\mathcal{L}_Q(\v,\w))$ corresponds to the weight
space of weight $\sum_{i \in I} (\w_i \omega_i - \v_i \alpha_i)$. In
\cite{Nak94}, Nakajima gave a similar realization of these
representations using a space of constructible functions on the
quiver varieties rather than their homology.  The irreducible
components of Nakajima quiver varieties enumerate a natural basis in
the representations of $\g_Q$.  These bases are given by the
fundamental classes of the irreducible components in the Borel-Moore
homology construction and by functions equal to one on an open dense
subset of an irreducible component (and equal to zero on an open
dense subset of all other irreducible components) in the
constructible function realization.


\section{Moduli spaces of representations of the Euclidean algebra}

Given that $\e(2)$ has wild representation type, it is prudent to
restrict one's attention to certain subclasses of modules and to
attempt a classification of the modules belonging to these classes.
One possible approach is to impose a restriction on the number of
generators of a representation (see \cite{Dou06a,Dou06b} for some
results in this direction and \cite{DP07} for other classes). In
this section we will examine the relationship between moduli spaces
of representations of the Euclidean algebra along with a set of
generating vectors and Nakajima quiver varieties.

Let $V$ be a finite-dimensional $\e(2)$-module.  For $u_1, u_2,
\dots, u_n \in V$, we denote by $\left<u_1, \dots, u_n \right>$ the
submodule of $V$ generated by $\{u_1, \dots, u_n\}$.  It is defined
to be the smallest submodule of $V$ containing all the $u_i$.  A
element $u \in V$ is called a \emph{weight vector} if it lies in
some weight space $V_k$ of $V$.  For a weight vector $u$, we let
$\wt u = k$ where $u \in V_k$.  We say that $\{u_1, \dots, u_n\}$ is
a set of \emph{generators} of $V$ if each $u_i$ is a weight vector
and $\left<u_1, \dots, u_n\right> = V$.  For $\v \in (\Z_{\ge
0})^\Z$, we let $|\v| = \sum_{k \in \Z} v_k$.

\begin{defin}
For $\v, \w \in (\Z_{\ge 0})^\Z$, let $E(\v,\w)$ be the set of all
\[
(V, (u^j_k)_{k \in \Z,\, 1 \le j \le \w_k})
\]
where $V$ is a finite-dimensional $\e(2)$-module with
$\mathbf{dim}\, V = \v$ and $\{u^j_k\}_{k \in \Z,\, 1 \le j \le
\w_k}$ is a set of generators of $V$ such that $\wt u^j_k = k$. We
say that two elements $(V, (u^j_k))$ and $(\tilde V, (\tilde
u^j_k))$ of $E(\v,\w)$ are \emph{equivalent} if there exists a
$\e(2)$-module isomorphism $\phi : V \to \tilde V$ such that
$\phi(u^j_k) = \tilde u^j_k$. We denote the set of equivalence
classes by $\E(\v,\w)$.
\end{defin}

\begin{theo} \label{thm:E=L}
There is a natural one-to-one correspondence between $\E(\v,\w)$ and
$\mathcal{L}_{Q_\infty}(\v,\w)$.
\end{theo}

\begin{proof}
Let $(V, (u^j_k)_{k \in \Z,\, 1 \le j \le \w_k}) \in E(\v,\w)$ and
let $V= \bigoplus V_k$ be the weight space decomposition of $V$.
Thus $V_k$ is isomorphic to $\C^{\v_k}$ and we identity the two via
this isomorphism.  We then define a point $\varphi(V,(u^j_k)) =
(x,s) \in L_{Q_\infty}(\v,\w)$ by setting
\begin{gather*}
x_{h_i} = p_+|_{V_i},\quad x_{\bar h_i} = p_-|_{V_{i+1}}, \quad i
\in \Z,\\
s(w^j_k) = u^j_k,\quad k \in \Z,\quad 1 \le j \le \w_k,
\end{gather*}
where $\{w^j_k\}_{1 \le j \le \w_k}$ is the standard basis of
$\C^{\w_k}$ and the map $s$ is extended by linearity.  It follows
from the results of Section~\ref{sec:euc-preproj} that $x \in
\Lambda_{V^\v,Q}$ and so $(x,s) \in L_{Q_\infty}(\v,\w)$.
Furthermore, it follows from the fact that $(u^j_k)$ is a set of
generators, that $(x,s)$ is a stable point.  Thus $\varphi :
E(\v,\w) \to L_{Q_\infty}(\v,\w)^{\text{st}}$.  It is easily
verified that two elements $(V,(u^j_k))$ and $(\tilde V,(\tilde
u^j_k))$ are equivalent if and only if $\varphi(V,(u^j_k))$ and
$\varphi(\tilde V,(\tilde u^j_k))$ lie in the same $G_\v$-orbit.
Thus $\varphi$ induces a map $\varphi' : \mathcal{E}(\v,\w) \to
\mathcal{L}_{Q_\infty}(\v,\w)$ which is independent of the
isomorphism $V \cong \C^{\v_k}$ chosen in our construction.  It is
easily seen that $\varphi'$ is a bijection.
\end{proof}

As noted in Section~\ref{sec:nak-qv}, the irreducible components of
Nakajima quiver varieties can be identified with the irreducible
components of Lusztig quiver varieties that are not killed by the
stability condition.  In the language of $\e(2)$-modules, passing
from Lusztig quiver varieties to Nakajima quiver varieties amounts
to imposing the condition that the module be generated by a set of
$|\w|$ weight vectors with weights prescribed by $\w$.

A \emph{partition} is a sequence of non-increasing natural numbers
$\lambda = (\lambda_1,\lambda_2,\dots,\lambda_l)$.  The
corresponding \emph{Young diagram} is a collection of rows of square
boxes which are left justified, with $\lambda_i$ boxes in the $i$th
row, $1 \le i \le l$.   We will identify a partition and its Young
diagram and we denote by $\mathcal{Y}$ the set of all partitions (or
Young diagrams). If $b$ is a box in a Young diagram $\lambda$, we
write $x \in \lambda$ and we denote the box in the $i$th column and
$j$th row of $\lambda$ by $x_{i,j}$ (if such a box exists).  The
\emph{residue} of $x_{i,j} \in \lambda$ is defined to be $\res
x_{i,j} = i-j$.  For $\lambda \in \mathcal{Y}$ and $a \in \Z$,
define $\v^{\lambda,a} \in (\Z_{\ge 0})^\Z$ by setting
$\v^{\lambda,a}_{i+a}$ to be the number of boxes in $\lambda$ of
residue $i$.

\begin{prop}
For $\lambda \in \mathcal{Y}$, there exists a unique $\e(2)$-module
$V$ (up to isomorphism) with a single generator of weight $a \in \Z$
and $\mathbf{dim}\, V = \v^{\lambda,a}$. It is given by
\begin{align*}
V &= \Span_\C \{x\ |\ x \in \lambda\} \\
l (x_{i,j}) &= (a+\res{x_{i,j}})x_{i,j} = (a + i - j)x_{i,j} \\
p_+ (x_{i,j}) &= x_{i+1,j} \\
p_- (x_{i,j}) &= x_{i,j+1},
\end{align*}
where we set $x_{i,j} = 0$ if there is no box of $\lambda$ in the
$i$th column and $j$th row.

For $\v \in (\Z_{\ge 0})^\Z$ such that $\v \ne \v^{\lambda,a}$ for
all $\lambda \in \mathcal{Y}$ and $a \in \Z$, there are no
$\e(2)$-modules $V$ with a single generator and $\mathbf{dim}\, V =
\v$
\end{prop}

\begin{proof}
By tensoring with an appropriate $\chi_n$, we may assume that the
generator of our module has weight zero.  It is shown in
\cite[\S5.1]{FS03} that
\[
\dim_\C \mathcal{L}_{Q_\infty}(\v,\w^0) = \begin{cases} 1 & \text{if
} \v = \v^{\lambda,0},\, \lambda \in \mathcal{Y},\\
0 & \text{otherwise} \end{cases},
\]
where $\w^0_0 = 1$ and $\w^0_i = 0$ for $i \ne 0$ (the first case
can be deduced from the dimension formula in
Proposition~\ref{prop:NQV-dim}). It then follows from
Theorem~\ref{thm:E=L} that if $V$ is an $\e(2)$-module with a single
generator $v$ of weight zero, we must have $\mathbf{dim}\, V =
\v^{\lambda,0}$. Furthermore, up to isomorphism, there is only one
such pair $(V,v)$ and thus only one such module $V$.
\end{proof}

Thus $\e(2)$-modules with a single generator of a fixed weight are
determined completely by the dimensions of their weight spaces. This
was proven directly by Gruber and Henneberger in \cite{GH83}. As in
the proof of Corollary~\ref{cor:finite-dim5}, we see that our
knowledge of the precise relationship between quivers and the
Euclidean algebra allows us to use known results about quivers and
quiver varieties to simplify such proofs.

\begin{rem}
As explained at the end of Section~\ref{sec:nak-qv}, the Nakajima
quiver varieties $\mathcal{L}_Q(\v,\w)$ are closely connected to the
representation theory of $\g_Q$. Therefore, the relationship noted
in Remark~\ref{rem:LQV-rep-e2} between the representation theory of
the Euclidean group and the Lie algebra $\mathfrak{sl}_\infty$ (or
the Lie groups $SL(n)$) is emphasized further by the above results.
Namely, the moduli space of representations of the Euclidean group
along with a set of generators is closely related to the
representation theory of $\mathfrak{sl}_\infty$ and the Lie groups
$SL(n)$.
\end{rem}

\begin{rem}
Although Theorem~\ref{thm:Euc-rep-type} tells us that the Euclidean
group has wild representation type, the results of this section
produce a method of approaching the unwieldy problem of classifying
its representations. Namely, if we fix the cardinality and weights
of a generating set, the resulting moduli space of representations
(along with a set of generators) is enumerated by a countable number
of finite-dimensional varieties, one variety for the representations
of each graded dimension.
\end{rem}


\section{Further directions}

The ideas presented in this paper open up some possible avenues of
further investigation.  We present here two of these.

Consider the Euclidean algebra over a field $k$ of characteristic
$p$ instead of over the complex numbers.  This algebra is still
spanned by $\{p_+,p_-,l\}$ with commutation relations
$\eqref{eq:e2-relations}$ but the weights of representations are
elements of $\Z/p\Z$ (if we restrict our attention to ``integral''
weights as usual) instead of $\Z$.  One can then show that this
category of representations is equivalent to the category of
representations of the preprojective algebra of the quiver of affine
type $\hat A_{p-1}$.  In this case, the representations with one
generator are, in general, more complicated than in the complex
case.  We refer the reader to \cite{FS03} for an analysis of the
corresponding quiver varieties.  There a graphical depiction of the
irreducible components of these varieties is developed.  These
quiver varieties are related to moduli spaces of solutions of
anti-self-dual Yang-Mills equations and Hilbert schemes of points in
$\C^2$ and it would be interesting to further examine the
relationship between these spaces and the Euclidean algebra.

In \cite{KS97} and \cite{S02}, Kashiwara and Saito defined a crystal
structure on the sets of irreducible components of Lusztig and
Nakajima quiver varieties. Using this structure, each irreducible
component can be identified with a sequence of crystal operators
acting on the highest weight element of the crystal. Under the
identification of quiver varieties with (moduli spaces of)
$\e(2)$-modules, these sequences correspond to the Jordan-H\"older
decomposition of $\e(2)$-modules.  It could be fruitful to further
examine the implications of this correspondence.


\bibliographystyle{abbrv}
\bibliography{biblist}

\end{document}